\documentclass{article}

\usepackage{latexsym}
\usepackage[centertags]{amsmath}
\usepackage{amsfonts}
\usepackage{amssymb}
\usepackage{amsthm}
\usepackage{newlfont}

\newcommand{\Complex}{\mathbb C}

\newcommand{\Int}{\mathbb Z}
\newcommand{\Nat}{\mathbb N}

\begin{document}

\begin{center}
{\LARGE \bf On the existence of polynomial first integrals  of
quadratic homogeneous systems of ordinary differential
equations}\\

 \vspace{0.3cm}

\ {\Large \bf  Tsygvintsev Alexei}

\end{center}

\begin{abstract}

We consider systems of ordinary differential equations with
quadratic homogeneous right hand side. We give a new simple proof
of  a result already obtained in [8,10] which gives the necessary
conditions for the existence of polynomial first integrals. The
necessary conditions for the existence of a polynomial symmetry
field are given. It is proved that an arbitrary homogeneous first
integral of a given degree is a linear combination of a fixed set
of polynomials.
\end{abstract}

\section{Introduction}

In the present paper we study the system of ordinary differential
equations with quadratic homogeneous right hand side

$$ \dot x_i=f_i(x_1,\ldots,x_n), \quad
f_i=\sum^n_{j,k=1}a_{ijk}x_j x_k , \quad a_{ijk} \in {\Complex},
\quad i=1,\ldots,n. \eqno (1.1) $$

Systems of such a form arise in many problems of classical mechanics:
Euler--Poincar\'e equations on Lie algebras, the Lotka--Volterra
systems, etc.

The main concern of this paper is to find the values of the parameters
$a_{ijk}$ for which equations (1.1) have first integrals.

In the paper [2] the necessary conditions are found for the
existence of polynomial first integrals of the system $$ \dot
x_i=V_i(x_1,\ldots,x_n), \quad i=1,\ldots,n, \eqno (1.2) $$ where
$V_i\in {\Complex}[x_1,\ldots,x_n]$ are homogeneous polynomials of
weighted degree $s\in \Nat$. In the case $s=2$ we obtain equations
(1.1).

The method given in [2] is based on ideas of Darboux [1,6,7] who
used a special type of particular solutions of the system (1.2) $$
 x_i(t)=d_i\phi(t),\quad i=1,\ldots,n,
$$
where $\phi(t)$ satisfies the differential equation
$\dot\phi=\lambda\phi^s$,
$\lambda$ is an arbitrary number and $d=(d_1,\ldots,d_n)^T\ne 0$ is
a solution of the following algebraic system
$$
V_i(d)=\lambda d_i,\quad i=1,\ldots,n.
$$

In this paper we generalize this
method.

 It was shown in [8,10] that the weighted degree of a polynomial first
integral of the system (1.1) is a certain integer linear
combination of Kovalevskaya exponents (see [9]). In Section 2 we
give a new simple proof of this result. In Section 3 a similar
theorem for polynomial symmetry fields is proved. As an example,
we consider the well known Halphen equations. Section 4 contains
our main result. We present so called { \it base functions } and
prove that every homogeneous polynomial first integral of a fixed
degree is a certain linear combination of the corresponding base
functions. In Section 5 we give an application of previous results
to planar homogeneous quadratic systems where necessary and
sufficient conditions for the existence of polynomial first
integrals in terms of Kovalevskaya exponents  are found.

\section{The existence of polynomial first integral. Necessary
conditions}

Following the paper [2], we consider  solution
$C=(c_1,\dots,c_n)^T\ne(0,\dots,0)^T$ of algebraic equations $$
f_i(c_1,\dots,c_n)+c_i=0,\quad i=1,\dots,n.\eqno (2.1) $$ Define
the {\it Kovalevskaya matrix} $K$ [3] $$ K_{ij}=\frac{\partial
f_i}{\partial x_j}(C)+\delta_{ij},\quad i,j=1,\dots,n. $$ where
$\delta_{ij}$ is the Kronecker symbol. Let us assume that $K$ can
be transformed to diagonal form $$
K=\mbox{diag}(\rho_1,\dots,\rho_n). $$
 The eigenvalues $\rho_1,\dots,\rho_n$ are called {\it Kovalevskaya
exponents}.
\newtheorem{lemma}{Lemma}
\begin{lemma}{\bf ([3])}
Vector C is an eigenvector of the matrix K with eigenvalue $\rho_1=-
1$.
\end{lemma}
Consider the following linear differential operators
$$
\begin{array}{ll}
D_+=\sum\limits^n_{i=1}f_i\frac{\displaystyle \partial}{\displaystyle
\partial x_i},
&\quad D_0=\sum\limits^n_{i,j=1}K_{ij}x_j
\frac{\displaystyle\partial}{\displaystyle \partial x_i},\\
U=\sum\limits^n_{i=1}x_i
\frac{\displaystyle\partial}{\displaystyle \partial x_i},&\quad
D_-=\sum\limits^n_{i=1}c_i
\frac{\displaystyle\partial}{\displaystyle \partial x_i},
\end{array} \eqno (2.2)
$$
which satisfy relations
$$
[D_-,D_+]=D_0-U,\quad [D_0,D_-]=D_-,\eqno (2.3)
$$
where $[A,B]=AB-BA$.
\newtheorem{theorem}[lemma]{Theorem}{\bf [8,10]}

\begin{theorem}
Suppose that the system (1.1) possesses a homogeneous polynomial
first integral $F_M$ of degree M, and
$\rho_1=-1,\rho_2,\dots,\rho_n$ are Kovalevskaya exponents. Then
there exists a set of non-negative integers $k_2,\dots,k_n$ such
that $$ k_2\rho_2+\cdots+k_n\rho_n=M,\quad k_2+\cdots+k_n\le M.
\eqno (2.4) $$
\end{theorem}
{\it Proof.} By definition of a first integral $D_+F_M=0$.
Considering identities $$ D^l_-(D_+F_M)=0,\quad \mbox{for}\quad
l\in \Nat, \eqno (2.5) $$ we obtain the following set of
polynomials $$ F_M,F_{M-1},\dots,F_{\rho+1},F_\rho, $$ defined by
the recursive relations $$ D_-F_{i+1}=(M-i)F_i,\quad
i=\rho,\dots,M-1, $$ where the number $1\le \rho \le M$ is
determined by the condition $D_- F_\rho=0$. Using (2.3), (2.5) we
deduce the following chain of equations $$
\begin{array}{l}
D_0F_M=MF_M-D_+F_{M-1},\\ D_0F_{M-1}=MF_{M-1}-D_+F_{M-2},\\
\dots\\ D_0F_{\rho+1}=MF_{\rho+1}-D_+F_\rho,\\ D_0F_\rho=MF_\rho.
\hspace{10.2cm} (2.6)
\end{array}
$$

Let $J_1$, \dots , $J_n$ -- linearly independent eigenvectors of
the Kovalevskaya matrix $K$ corresponding to the eigenvalues
$\rho_1=-1$, $\rho_2$, \dots , $\rho_n$. According to Lemma 1 we
can always put $J_1=C$.

Consider the linear change of variables $$
x_i=\sum^n_{j=1}L_{ij}y_j,\quad i=1,\dots,n, \eqno (2.7) $$ where
$L=(L_{ij})$ is a nonsingular matrix defined by $$
L=(C,J_2,\dots,J_n),$$ then obviously $$
L^{-1}KL=\mbox{diag}(-1,\rho_2,\dots,\rho_n). $$

With help of (2.7) and Lemma 1 one finds the following expressions
for the operators  $D_0$, $D_-$ in the new variables $$
D_0=\sum^n_{i=1}\rho_iy_i\frac{\partial}{\partial y_i}, \quad
D_-=-\frac{\displaystyle
\partial}{\displaystyle
\partial y_1},
$$ and the equation (2.6) becomes $$
\left(\rho_2y_2\frac{\partial}{\partial y_2}+\cdots+\rho_ny_n
\frac{\partial}{\partial y_n}\right)F_\rho=MF_\rho. \eqno (2.8) $$
We can write the polynomial $F_\rho$ as follows $$
F_\rho=\sum_{|k|=\rho}A_{k_2\dots k_n}y^{k_2}_2\cdots y^{k_n}_n,
\quad |k|=k_2+\cdots+k_n,\quad k_i \in \Int_+. \eqno (2.9) $$
Substituting (2.9) into (2.8), one obtains the following linear
system $$ (k_2\rho_2+\cdots+k_n\rho_n)A_{k_2\dots k_n}=MA_{k_2
\dots k_n}, \quad \mbox{for} \quad |k|=\rho. \eqno (2.10) $$
Taking into account that $F_\rho$ is not zero identically, we
conclude that there exists at least one nonzero set
$k_2,\dots,k_n$, $|k| \le M $ such that $$
k_2\rho_2+\cdots+k_n\rho_n=M. \eqno (2.11) $$ This relation
implies (2.4). Q.E.D.

\newtheorem{cor}[lemma]{Corollary}

\noindent {\it Remark.}  Theorem 2 does not impose any
restrictions on $\mbox{grad}(F_M)$ calculated at the point C.
Thus, it generalizes the theorem of Yoshida ([3], p.572), who used
essentially the condition $\mbox{grad}(F_M) \ne 0$.

\begin{cor} {\rm
The Halphen equations
$$
\begin{array}{l}
\dot x_1=x_3x_2-x_1x_3-x_1x_2,\\
\dot x_2=x_1x_3-x_2x_1-x_2x_3,\\
\dot x_3=x_2x_1-x_3x_2-x_3x_1,
\end{array} \eqno (2.12)
$$
 admit no polynomial first integrals.}
\end{cor}
Indeed, the system (2.12) has Kovalevskaya exponents
$\rho_1=\rho_2=\rho_3=-1$. It is easy to verify that conditions
(2.4) are not fulfilled for any positive integer $M$. Moreover, as
proved in [2], the system (2.12) has no rational first integrals.

\section{Existence of polynomial symmetry fields. Necessary
conditions}

The first integrals are the simplest tensor invariants of the system
(1.1). In [4] Kozlov
considered tensor invariants of weight-homogeneous differential
equations which include
the system (1.1). In particular, he found necessary conditions for
the existence
of symmetry fields. Below we propose a generalization of his result.

Recall that the linear operator
$
W=\sum\limits^n_{i=1}w_i(x_1,\dots,x_n) \frac{\displaystyle
\partial}{\displaystyle \partial x_i},
$
is called the {\it symmetry field} of (1.1), if $[W,D_+]=0$,where
$D_+$ is defined by (2.2). If $w_1,\dots,w_n$ are homogeneous
functions  of degree $M+1$ then the degree of $W$ is $M$ ([4]).
\begin{theorem}
Suppose that the system (1.1) possesses a polynomial symmetry
field of degree $M$ and \\ $\rho_1=-1,\rho_2,\dots,\rho_n$ are
Kovalevskaya exponents. Then there exist non-negative integers
$k_2,\dots,k_n$, $|k|\le M+1$ such that at least one of the
following equalities holds $$
k_2\rho_2+\cdots+k_n\rho_n=M+\rho_i,\quad i=1,\dots,n. \eqno (3.1)
$$
\end{theorem}
{\it Proof.}
Let $W_M$ be a polynomial symmetry field of degree $M$.
Substituting in the proof of Theorem 2 operators $D_+$, $D_0$, $D_-$
with their commutators $[D_+,{}]$, $[D_0,{}]$, $[D_-,{}]$
respectively we repeat the same arguments.

\begin{cor}
{\rm The Halphen equations (2.12) admit no polynomial symmetry
fields.}
\end{cor}
Using (3.1) we obtain that $M$ may be equal to $-1$, $0$, $1$ only.
It is easy to check that (2.12) does not have symmetry fields
of such degrees.

\section{Base functions}

After the change of variables (2.7) the system (1.1) takes the form
$$
\begin{array}{l}
\dot y_1=-y^2_1+\varphi_1(y_2,\dots,y_n)\\ \dot
y_i=(\rho_i-1)y_1y_i+\varphi_i(y_2,\dots,y_n), \quad i=2,\dots,n,
\end{array} \eqno (4.1)
$$ where $\varphi_i$ are quadratic homogeneous polynomials in the
variables $y_2,\dots,y_n$.

According to (2.2) define operators $D_+,D_0,D_-$.

A homogeneous polynomial $P_M(y_1,\dots,y_n)$ of degree $M$
satisfying the condition $$ D_-(D_+P_M)=0, \eqno(4.2) $$ is called
{\it base function} of the system (4.1). In others words, the
function $D_+P_M$ does not depend on $y_1$. It is clear that base
functions of degree $M$ form a linear space $L_M$ over field
$\Complex$.
\begin{lemma}
If the system (4.1) has a homogeneous polynomial first integral $F_M$
of degree $M$, then $F_M\in L_M$.
\end{lemma}
Indeed, by definition, we have $D_+F_M=0$, hence, in view of (4.2),
$F_M\in L_M$.

Let $J(M)=\{z\in {\Int^{n-1}_+}\mid z_2\rho_2+\cdots+z_n\rho_n=M,
|z|\leq M\}$ be the set of integer-valued vectors
$z=(z_2,\dots,z_n)^T$ for which the condition (2.11) is fulfilled.
Put $m=|J(M)|$ and suppose $J(M)\ne \emptyset$.
\begin{theorem}
The  dimension $d$ of $L_M$ satisfies the  condition $1\le d \le m$.
\end{theorem}
{\it Proof.} Let us assume the set J(M) contains vectors
$z^{(1)},\dots,z^{(m)}$ which are ordered by the norm
\\ $|z|=z_2+\cdots+z_n$ $$ |z^{(1)}|\le \cdots \le |z^{(m)}|. $$
Define the vector $\rho=(\rho_2,\dots,\rho_n)^T$ and  put
$(\rho,z)=\rho_2z_2+\cdots+\rho_nz_n$,
 $|z^{(i)}|=n_i$, $i=1,\dots,m$.

Following the proof of Theorem 2, for each $i=1,\dots,m$ consider
the system of linear partial differential equations $$
\begin{array}{l}
D_0P_{i,n_i}=MP_{i,n_i},\\
D_0P_{i,n_i+1}=MP_{i,n_i+1}-D_+P_{i,n_i},\\
\cdots\\
D_0P_{i,M-1}=MP_{i,M-1}-D_+P_{i,M-2},\\
D_0P_{i,M}=MP_{i,M}-D_+P_{i,M-1},\\
\end{array} \eqno (4.3)
$$ $$ D_-P_{i,l+1}=(M-l)P_{i,l},\quad l=n_i,\dots,M-1,\eqno (4.4)
$$
 which defines polynomials $P_{i,n_i},\dots,P_{i,M}$ recurrently.

It follows from $(z^{(i)},\rho)=M$ that the first equation in (4.3)
has the
particular solution
$
P_{i,n_i}=y^{z_2^{(i)}}_2\cdots y^{z_n^{(i)}}_n.
$

Equations (4.3), (4.4) define certain base function $P_{i,M}$.
Indeed, according to (4.4), we have $$ P_{i,M-1}=D_-P_{i,M}. \eqno
(4.5) $$ Substituting (4.5) into the last equation in (4.3), and
using relations (2.3) we find $$
D_0P_{i,M}=MP_{i,M}-D_+D_-P_{i,M}=MP_{i,M}-(D_-D_+-D_0+U)P_{i,M}.
$$ Hence $D_-(D_+P_{i,M})=0$.

Now consider the problem on the existence of a solution of (4.3),
(4.4) in  form of homogeneous polynomials
$P_{i,n_i},\dots,P_{i,M}$. Fix certain $i=1,\dots,m$ and put
$a_i=M-n_i$. Using the relations (4.4) we can write $$
\begin{array}{l}
P_{i,n_i}=I_{i,n_i},\\ P_{i,n_i+p}=\sum\limits^p_{j=0}{p-j\choose
a_i-j}y^{p- j}_1I_{i,n_i+j},\quad p=1,\dots,a_i,
\end{array}  \eqno (4.6)
$$ where $I_{i,k}(y_2,\dots,y_n)$ are certain homogeneous
polynomials of degrees $k=n_i,\dots,M$. Notice that $I_{i,k}$ does
not depend on $y_1$.

Differential operators $D_+$, $D_0$ can be represented in the form
$$
\begin{array}{l}
D_+=(-y^2_1+\varphi_1)\frac{\displaystyle \partial}{\displaystyle
\partial y_1}+y_1(A_0-\tilde U)+A_+,\\
D_0=-y_1\frac{\displaystyle \partial}{\displaystyle \partial y_1}+A_0,
\end{array} \eqno (4.7)
$$
where
$$
A_+=\sum^n_{k=2}\varphi_k\frac{\partial}{\partial y_k},\quad
A_0=\sum^n_{k=2}\rho_ky_k
\frac{\partial}{\partial y_k},\quad
\tilde U=\sum^n_{k=2}y_k\frac{\partial}{\partial y_k}.
\eqno (4.8)
$$
Using (4.3), (4.6), (4.7) one deduces the following equations for
determination
of $I$
$$
\begin{array}{l}
A_0I_{i,n_i}=MI_{i,n_i},\\
A_0I_{i,n_i+1}=MI_{i,n_i+1}-A_+I_{i,n_i},\\
A_0I_{i,n_i+2}=MI_{i,n_i+2}-a_i\varphi_1I_{i,n_i}-A_+I_{i,n_i+1},\\
A_0I_{i,n_i+3}=MI_{i,n_i+3}-(a_i-1)\varphi_1I_{i,n_i+1}-
A_+I_{i,n_i+2},\\
\cdots\\
A_0I_{i,M}=MI_{i,M}-2\varphi_1I_{i,M-2}-A_+I_{i,M-1}.
\end{array} \eqno (4.9)
$$ We can write each equation of (4.9) as follows $$
A_0X_l=MX_l+Y_l, \eqno (4.10) $$ where $X_l,Y_l$ are homogeneous
polynomials of weighted degree $l=n_i,\dots,M$. Let us assume $$
X_l=\sum\limits_{|i|=l}c_{i_2\dots i_n}y^{i_2}_2\cdots y^{i_n}_n,
\quad Y_l=\sum\limits_{|i|=l}d_{i_2\dots i_n}y^{i_2}_2\cdots
y^{i_n}_n,\quad |i|=i_2+\cdots+i_n, \eqno (4.11) $$ where
$c_{i_2\dots i_n}$, $d_{i_2\dots i_n}$ are constant parameters.
Then  substituting (4.11) into (4.10), we obtain the following
linear system with respect to $c_{i_2\dots i_n}$ $$
(i_2\rho_2+\cdots+i_n\rho_n-M)c_{i_2\dots i_n}=d_{i_2\dots i_n},
\eqno (4.12) $$ for $i_2,\dots,i_n=0,1,\dots$, $|i|=l$.

Suppose there exists a set $k_2,\dots,k_n$ for which the following
conditions are fulfilled $$ (k_2,\dots,k_n)^T\in J(M),\quad
d_{k_2\dots k_n}\ne 0,\quad |k|=l, \eqno (4.13) $$ then the
solution $I_{i,n_i},\dots,I_{i,M}$ does not exist. In this case we
put $P_{i,M}=0$.

If the conditions (4.13) are not satisfied, we obtain the base
function $$ P_{i,M}=\sum^{a_i}_{j=0}y^{a_j-j}_1I_{i,n_i+j}.\eqno
(4.14) $$ It is easy to show that polynomials
$\{P_{i,M}\}^{i=1}_{i=m}$ are linearly independent over the field
$\Complex$.

Taking into account that $n_1\le\cdots\le n_m$ and using (4.13),
we see that in case $i=m$ we always can determine the base
function $P_{i,M}$. Therefore, under the assumption
$J(M)\ne\emptyset$, the space $L_M$ always contains a nonzero
function. Q.E.D.

\begin{cor} {\rm If at least one resonance condition of the form $$
(z,\rho)=M,\quad |z|\le M,\quad z\in {\Int^{n-1}_+}, $$ is
fulfilled, then there exists a base function of degree $M$.}
\end{cor}

\section{ Polynomial first integrals in the case of quadratic plane vector field.}
The first classification of integral curves of two-dimensional quadratic homogeneous systems
can be found in the paper of  Lyagina [5] and later was completed by numerous authors.

In this section we apply the previous results to this problem to
illustrate the method of basis functions.

Consider the system $$
\begin{array}{l}
\dot x_1=a_1x^2_1+b_1x_1x_2+d_2x^2_2,\\
\dot x_2=a_2x^2_2+b_2x_1x_2+d_1x^2_1,
\end{array} \eqno (5.1)
$$ where $a_i,b_i,d_i$ are constant parameters.

Let $c^{(1)}=(c^{(1)}_1,c^{(1)}_2)^T$,
$c^{(2)}=(c^{(2)}_1,c^{(2)}_2)^T$ be any two linearly independent
solutions of the algebraic system (2.1). The exceptional cases
when the system (2.1) has only one or admit no solutions are
excluded for the discussion below.

Assume that Kovalevskaya exponents corresponding to $c^{(1)}$,
$c^{(2)}$ are $$ \Re_1=(-1,\rho_1)^T,\quad
\Re_2=(-1,\rho_2)^T.\eqno (5.2) $$
\begin{lemma}
The system (5.1) has a homogeneous polynomial first integral of
degree M if and only if there exists integer $k=1,\dots,M-1$ such
that $\rho_1=M/k$ and $\rho_2$ is one of the following numbers $$
\frac{M}{M-k},\frac{M}{M-k-1},\dots,\frac{M}{2},M. $$
\end{lemma}
{\it Proof.}
Consider the following change of coordinates
$$
\left(
\begin{array}{l}
x_1\\
x_2
\end{array}
\right)
=
\left(
\begin{array}{ll}
c^{(1)}_1&c^{(2)}_1 \\
c^{(1)}_2&c^{(2)}_2
\end{array}
\right)
\left(
\begin{array}{l}
p_1\\
p_2
\end{array}
\right)\eqno (5.3)
$$
which exists because of linear independence of vectors $c^{(1)}$,
$c^{(2)}$.
In coordinates $(p_1,p_2)$ the system (5.1) takes a more simple form
$$
\begin{array}{l}
\dot p_1=-p^2_1+(\rho_2-1)p_1p_2,\\
\dot p_2=-p^2_2+(\rho_1-1)p_1p_2.
\end{array} \eqno (5.4)
$$ It is easy to show that under the change (5.3), the  vectors
$c^{(1)}$, $c^{(2)}$ turn into $\tilde c^{(1)}=(1,0)^T$, $\tilde
c^{(2)}=(0,1)^T$ respectively. Obviously, the system (5.4) has the
same Kovalevskaya exponents (5.2). The matrix $K$, calculated for
$\tilde c^{(1)}$ is $$ K= \left(
\begin{array}{ll}
-1&\rho_2-1\\
0&\rho_1\\
\end{array}
\right)
$$

Under the assumption $\rho_1\ne -1$, we can reduce $K$ to a diagonal
form using the following
change of coordinates
$$
\left(
\begin{array}{l}
p_1\\
p_2
\end{array}
\right)=L
\left(
\begin{array}{l}
y_1\\
y_2
\end{array}\right)
$$
with the constant matrix $L$
$$
L=
\left(
\begin{array}{ll}
1&\rho_2-1\\
0&\rho_1-1
\end{array}
\right)
$$
The case $\rho_1=-1$ will be considered below.

Finally, equations (5.4) take the form (4.1)
$$
\begin{array}{l}
\dot y_1=-y_1^2+\varphi_1,\\
\dot y_2=(\rho_1-1)y_1y_2+\varphi_2,
\end{array} \eqno (5.5)
$$
where
$$
\begin{array}{l}
\varphi_1=ay^2_2, \quad
\varphi_2=by^2_2,\\
a=(\rho_2-1)(\rho_1+\rho_2),\quad b=(\rho_1-1)(\rho_2-1)-\rho_1-1.
\end{array}
$$
For the operators (4.8) we get
$$
A_+=\varphi_2\frac{\partial}{\partial y_2},\quad
A_0=\rho_1y_2\frac{\partial}{\partial y_2},\quad
\tilde U=y_2\frac{\partial}{\partial y_2}.
$$
Let $F_M$ be a polynomial first integral of (5.5) of degree
$M$.

According to Theorem 2, there exists an integer $k=1,\dots,M-1$
such that $$ k\rho_1=M. \eqno (5.6) $$ We exclude the case $k=M
(\rho_1=1)$, since if $\rho_1=\pm 1$, then the system (5.1) has no
polynomial first integrals. This can be shown directly using
equations (5.4), (5.5).

Next, calculate the base function $P_M$ corresponding to the
resonance condition (5.6). Consider the equations (4.9). It is
obvious, that polynomials $I_{1,k},\dots,I_{1,M}$  can be
represented in the following form $$
I_{1,k+i}=\alpha_iy^{k+i}_2,\quad i=0,\dots,M-k, \eqno (5.7) $$
where $\alpha_0,\dots,\alpha_{M-k}$ are constant parameters.

Substituting (5.7) into (4.9) we obtain
$$
\begin{array}{l}
\alpha_0=1,\\ \alpha_1=\frac{\displaystyle bk}{\displaystyle
\rho_1(k+1)-M},\\ \alpha_i=\frac{\displaystyle
a(M-k-i+2)\alpha_{i-1}+b(k+i-1) \alpha_{i-2}}{\displaystyle
\rho_1(k+i)-M},\quad i=2,\dots,M-k.
\end{array} \eqno (5.8)
$$
According to (4.14), we get the following expression for the base
function $P_M$
$$
P_M=\sum^{M-k}_{j=0}\alpha_j y^{M-k-j}_1y_2^{k+j}. \eqno (5.9)
$$
By definition of the base function it is clear that
$$
D_+P_M=\delta y_2^{M+1}, \eqno (5.10)
$$
where
$$
\delta=a\alpha_{M-1}+bM\alpha_M. \eqno (5.11)
$$
Thus, the linear space $L_M$ contains only one polynomial $P_M$.
Hence, taking into account Lemma 6, $F_M=\mbox{const}\cdot P_M$.

Using (5.10), we conclude that $P_M$ is a first integral if and
only if $\delta =0$. In view of (5.6), (5.8), (5.11) and the above
condition, we arrive at Lemma 9.

\begin{theorem}

The system (5.1) possesses a homogeneous polynomial first integral
of degree $M$ if and only if the following conditions are
fulfilled \\

\noindent a) $\rho_i$, $i=1,2$ -- are positive rational numbers,
\\

\noindent b) $\rho_1^{-1}+\rho_2^{-1}\leq 1$, \\

\noindent c) $\displaystyle \frac{M}{\rho_i}\in \Nat.$

\end{theorem}

This is an obvious consequence of Lemma 9.

As an example consider the following system $$
\begin{array}{l}
\dot x_1=x_1^2-9x^2_2,\\
\dot x_2=-3x^2_1-8x_1x_2+3x_2^2.
\end{array} \eqno (5.12)
$$ The vectors $c^{(1)}$, $c^{(2)}$ have the form $$
c^{(1)}=(1/8,-1/8)^T,\quad c^{(2)}=(1/8,1/8)^T. $$ Calculating the
corresponding  Kovalevskaya exponents (5.2) one obtains $$
\Re_1=(-1,3)^T,\quad \Re_2=(-1,3/2)^T. $$

We have $\rho_1=3$, $\rho_2=3/2$, $\rho_1^{-1}+\rho_2^{-1}=1$. So,
the conditions a), b) of Theorem 10 are fulfilled. By the
condition c) one gets $M=3l$, $l\in \Nat$.
 Thus, the equations (5.12) possess a cubic first integral $F_3$.
Using formulas (5.8), (5.9), we obtain $$
F_3=x_1^3+x_1^2x_2-x_1x^2_2-x_2^3. $$

\vspace{1cm}
\noindent {\bf Acknowledgments}

We thank J.-M. Strelcyn for his useful remarks, V. Kozlov, D.
Treshev, H. Yoshida, K. Emelyanov, Yu. Fedorov and L. Gavrilov for
the interest to the paper.

\vspace{0.3cm}

 \noindent{\small \bf
 Section de Mathematiques,\\
 Universit\'e de Gen\`eve\\
 2-4, rue du Lievre,\\
 CH-1211, Case postale 240, Suisse \\
Tel l.: +41 22 309 14 03 \\
 Fax: +41 22 309 14 09 \\
E--mail: Alexei.Tsygvintsev@math.unige.ch}

\end{document}